\documentclass[leqno]{bit2e}

\usepackage{amssymb}
\usepackage{eqnarray}

\def\bitinfo{  }

\usepackage{epsfig} 

\newcommand{\beq}[2]{\begin{equation}
              \label{#1}
                {#2}
           \end{equation}}
\newcommand{\bea}{\begin{array}}
\newcommand{\ea}{\end{array}}
\newcommand{\beea}{\begin{eqnarray}}
\newcommand{\eea}{\end{eqnarray}}
\def\q{\quad}
\def\mr#1{$(\ref{#1})$}
\newcommand{\Matrix}[1]{\left[ \matrix{#1} \right]}
\newcommand{\Ac}{\mathcal A}
\newcommand{\Bc}{\mathcal B}
\newcommand{\Dc}{\mathcal D}
\newcommand{\Ec}{\mathcal E}

\newcommand{\Ic}{\mathcal I}
\newcommand{\Kc}{\mathcal K}
\newcommand{\Lc}{\mathcal L}
\newcommand{\Mc}{\mathcal M}
\newcommand{\Oc}{\mathcal O}
\newcommand{\Rc}{\mathcal R}
\newcommand{\Uc}{\mathcal U}
\newcommand{\Vc}{\mathcal V}
\newcommand{\real}{\mathbb{R}}
\newcommand{\complex}{\mathbb{C}}
\newcommand{\sinv}{\frac{1}{s}}
\newcommand{\szinv}{\frac{1}{s_0}}

\newcommand{\Spa}{\mathop{\rm span}\nolimits}
\newcommand{\ignore}[1]{}
\title{{\Large KRYLOV SUBSPACES ASSOCIATED WITH\\[4pt]
               HIGHER-ORDER LINEAR DYNAMICAL SYSTEMS}
}
\author{\small ROLAND W. FREUND}
\date{}

\begin{document}
\maketitle

\baselineskip=0.9\normalbaselineskip
\vspace{-3pt}

\begin{center}
{\footnotesize\em 
Department of Mathematics,
University of California at Davis,
One Shields Avenue,\\
Davis, California 95616, U.S.A.\\
email: freund\symbol{'100}math.ucdavis.edu}
\end{center}

\begin{abstract}
A standard approach to model reduction of large-scale higher-order
linear dynamical systems is to rewrite the system as an
equivalent first-order system and then employ Krylov-subspace
techniques for model reduction of first-order systems.
This paper presents some results about the structure
of the block-Krylov subspaces induced by the matrices of
such equivalent first-order formulations of higher-order
systems.
Two general classes of matrices, which exhibit
the key structures of the matrices of first-order 
formulations of higher-order systems, are introduced.
It is proved that for both classes, the block-Krylov subspaces induced by the matrices
in these classes can be viewed as multiple copies of certain subspaces of the 
state space of the original higher-order system.
\end{abstract}

\medskip
\begin{AMS}
{\small 65F30, 15A57, 65P99, 41A21.} 
\end{AMS}
  
\medskip    

\baselineskip=\normalbaselineskip

\section{Introduction.} \label{sec-introduction}

In recent years, Krylov-subspace methods, especially the Lanczos algorithm
and the Arnoldi process, have become popular tools for 
model reduction of large-scale time-invariant linear dynamical systems;
we refer the reader to the survey papers~\cite{Fre97b,Fre00a,Bai02,Fre03b},
and the references given there.
Krylov-subspace techniques can be applied directly only to
first-order linear dynamical systems.
However, there are important applications, for example in VLSI circuit
simulation \cite{Tut77,Fre04c,Fre05a}, structural dynamics \cite{Prz85,CraH88,SuC91},
and computational electromagnetics~\cite{WitSW03}, that lead to second-order,
or even general higher-order, linear dynamical systems.

The standard approach to employing Krylov-subspace
methods for model reduction of a second-order or 
higher-order system is to first rewrite the system as an
equivalent first-order system, and then apply Krylov-subspace
techniques for reduced-order modeling of first-order systems.
At first glance, there are two disadvantages of this standard approach.
First, the second-order or higher-oder structure is not preserved
by a straightforward application of Krylov-subspace methods to the
first-order formulation.
Second, the computational cost increases due to the fact that the
state-space dimension of the first-order formulation is $l$ times
the state-space dimension of the original $l$-th-order system.
A partial remedy of the first problem is to use certain 
structure-preserving projections, as described in the recent
papers \cite{VanV04,Fre04c,Van04,Fre05a,ChaGVV05}.
However, the structure-preserving property of these approaches comes at
the expense of reduced approximation quality of the resulting
models.
To address the second problem, at least for the special case of 
second-order systems, various authors have proposed to directly generate
basis vectors of certain subspaces of the state space of the second-order
system, rather than basis vectors of the Krylov subspaces of the first-order
formulation; see, e.g.,~\cite{BaiS03,Li04,RenZ97,SuC91,ZheRW97}.

The purpose of this paper is to shed some light on the second problem
and to present some results on the special structures of the block-Krylov
subspaces induced by the matrices of equivalent first-order formulations
of general higher-order time-invariant linear dynamical systems and
of certain systems of first-order integro-differential-algebraic equations.
More precisely, we introduce two classes of structured matrices, which
include the matrices of these first-order
formulations as special cases.
As our main results, we show that the block-Krylov subspaces
induced by the matrices in theses classes exhibit special structures.
Roughly speaking, for both classes, the associated structured
block-Krylov subspaces consist of multiple copies of certain
subspaces of the state space of the original higher-order system.

The remainder of the paper is organized as follows.
In Section~\ref{sec-krylov}, we briefly review the notion of block-Krylov subspaces.
In Section~\ref{sec-main}, we introduce two classes of matrices,
and we state our main results about the special structures
of the block-Krylov subspaces associated with these two classes.
In Section~\ref{sec-proofs}, we present proofs of these main results.
In Section~\ref{sec-higher}, we consider higher-order linear dynamical systems,
and we show how certain model-reduction approaches lead to matrices that
are special instances of the first class of matrices introduced 
in~Section~\ref{sec-main}.
In Section~\ref{sec-second}, we study systems of first-order
integro-differential-algebraic equations,
and we show how model reduction leads to matrices that
are special instances of the second class of matrices introduced 
in~Section~\ref{sec-main}.
Finally, in Section~\ref{sec-concluding}, we make some concluding remarks.

Throughout this paper the following notation is used.
The set of real and complex numbers is denoted by $\real$ and
$\complex$, respectively.
Unless stated otherwise, all vectors and matrices are allowed to have real
or complex entries.
For a matrix $A = \Matrix{a_{jk}}\in \complex^{m\times n}$, we denote by
$A^H := \Matrix{\overline{a_{kj}}} \in \complex^{n \times m}$
its conjugate transpose.
For any two matrices 
$A = \Matrix{a_{jk}}\in \complex^{m\times n}$ and $B\in \complex^{p\times q}$,
\[
A \otimes B := \Matrix{a_{jk} B}\in \complex^{mp \times nq}
\]
is the Kronecker product~\cite{Gra81,Van00} of $A$ and $B$.
The $n\times n$ identity matrix is denoted by $I_n$ and the
zero matrix by $0$.
If the dimension of $I_n$ is apparent from the context,
we drop the index and simply use $I$.
The actual dimension of $0$ will always be
apparent from the context.

\section{Block-Krylov subspaces.} \label{sec-krylov} 

We use the notion of block-Krylov subspaces that was introduced
in~\cite{AliBFH00} in connection with a band Lanczos process
for multiple starting vectors.
In this section, we briefly review the definition of block-Krylov 
subspaces from~\cite{AliBFH00}.

In the following, let 
\beq{KA0}{
\Mc \in \complex^{N\times N}\q \mbox{and}\q
\Rc =\Matrix{r_1 & r_2 & \cdots & r_m}\in \complex^{N\times m}
}
be given matrices.
The $N\times mN$ matrix
\beq{KAA}{
\Matrix{\Rc & \Mc\Rc & \Mc^2 \Rc & \cdots & \Mc^{N-1} \Rc}
}
is called the \textit{block-Krylov matrix} induced by $\Mc$ and $\Rc$.

\subsection{The case of exact deflation.} \label{ssec-exact}

Let $N_0 \, (\leq N)$ denote the rank of the block-Krylov 
matrix~\mr{KAA}.
Hence only $N_0$ of the $mN$ columns of~\mr{KAA}
are linearly independent.
Such a set of $N_0$ linearly independent columns
can be constructed by scanning the columns of~\mr{KAA} 
from left to right and deleting each column that is linearly
dependent on earlier columns.
This process of deleting linearly dependent columns is called
\textit{exact deflation}.
By the structure of the block-Krylov matrix~\mr{KAA}, a column $\Mc^{k-1} r_i$
being linearly dependent on earlier
columns implies that all columns $\Mc^j r_i$, $k\leq j \leq N-1$,
are also linearly dependent on earlier columns.
Consequently, applying exact deflation to~\mr{KAA} results 
in a matrix of the form
\beq{KAC}{
\Vc(\Mc,\Rc) :=
\Matrix{\Rc_1 & \Mc\Rc_2 & \Mc^2 \Rc_3 & \cdots & \Mc^{k_0-1} \Rc_{k_0}}
\in \complex^{N\times N_0}.
}
Here, for each $k=1,2,\ldots,k_0$,
$\Rc_k \in \complex^{N\times m_k}$ is a submatrix of 
$\Rc_{k-1} \in \complex^{N\times m_{k-1}}$, 
with $\Rc_k\not=\Rc_{k-1}$ if, and only if,
exact deflation occurs within the $k$-th Krylov block $\Mc^{k-1} \Rc$
in~\mr{KAA}.
(For $k=1$, we set $\Rc_0=\Rc$ and $m_0 = m$.)
For later use, we remark that 
\beq{KAD}{
\Rc_k = \Rc_{k-1} E_k,\q E_k \in \complex^{m_{k-1} \times m_k},\q m_k \leq m_{k-1},
}
where $E_k$ is the \textit{deflated identity matrix} obtained from $I_{m_{k-1}}$ by 
deleting those $m_{k-1}-m_k$ columns corresponding to exact deflation within 
the $k$-th Krylov block.

By construction, the matrix~\mr{KAC} has full column rank $N_0$.
For $1\leq n\leq N_0$, the \textit{$n$-th block-Krylov subspace}
(induced by $\Mc$ and $\Rc$), $\Kc_n(\Mc,\Rc)$, is defined as the
$n$-dimensional subspace of $\complex^N$ spanned by the first $n$
columns of~\mr{KAC}.
We say that  
\[
\Vc = \Matrix{v_1 & v_2 & \cdots & v_{N_0}} \in \complex^{N\times N_0}
\]
is a \textit{basis matrix} of the block-Krylov subspaces induced by $\Mc$ 
and $\Rc$ if
\[
\Kc_n(\Mc,\Rc) = \Spa \Matrix{v_1 & v_2 & \cdots & v_n}\q \mbox{for all}\q
n=1,2,\ldots,N_0.
\]
Note that the matrix $\Vc(\Mc,\Rc)$ defined in~\mr{KAC} is a particular
instance of a basis matrix.
Furthermore, any two basis matrices $\Vc_1$ and $\Vc_2$ of the block-Krylov 
subspaces induced by $\Mc$ and $\Rc$ are connected by a relation of the form
\beq{KAM}{
\Vc_1 = \Vc_2\, \Uc, 
}
where $\Uc$ is a nonsingular and upper triangular matrix.

Lanczos- and Arnoldi-type algorithms for the actual construction of 
basis matrices of block-Krylov subspaces can be found in~\cite{AliBFH00} 
and~\cite{Fre03b}.  

\subsection{Inexact deflation.}

In the above construction of block-Krylov
subspaces, we performed only exact deflation.
In an actual algorithm for constructing a basis matrix
of the block-Krylov subspaces induced by $\Mc$ 
and $\Rc$ in finite-precision arithmetic,
one also needs to delete vectors that are in some
sense ``almost'' linearly dependent on earlier vectors.
The deletion of such almost linearly dependent
vectors is called \textit{inexact deflation}. 
For example, the Lanczos- and Arnoldi-type algorithms
in~\cite{AliBFH00} and~\cite{Fre03b} have simple built-in procedures for 
both exact and inexact deflation.

It turns out that the construction of block Krylov
subspaces described in Subsection~\ref{ssec-exact} can be extended
to the more general case when exact and inexact deflations are performed.
The deflated matrix~\mr{KAC} is now obtained by deleting from
the block-Krylov matrix~\mr{KAA} those columns that are linearly or
almost linearly dependent on columns to their left.
In the general case, $N_0$ is now simply defined as the number of columns 
of the resulting deflated matrix~\mr{KAC}.
Note that $N_0$ is less than or equal to the rank of the block-Krylov 
matrix~\mr{KAA}, with equality only if no inexact deflation occurs.
Based on the deflated matrix~\mr{KAC}, block-Krylov subspaces and basis
matrices of these subspaces are defined in the same way as 
in Subsection~\ref{ssec-exact}.
However, note that the resulting block-Krylov subspaces are in general
different from the block-Krylov subspaces obtained with exact deflation
only.  

The main results of this paper, namely Theorems~\ref{theorem1} 
and~\ref{theorem2} below, hold true for the general case of 
exact and inexact deflations, provided that the matrices 
$\Rc_k$ in~\mr{KAC} still satisfy relations of the form~\mr{KAD}.
This is the case for the built-in deflation procedures of the
Lanczos- and Arnoldi-type algorithms 
in~\cite{AliBFH00} and~\cite{Fre03b}.
Thus, in the following, we always assume that the matrices 
$\Rc_k$ in~\mr{KAC} indeed satisfy relations of the form~\mr{KAD}.

\section{Main results.} \label{sec-main}
 
In this section, we introduce two classes of matrices $\Mc$ and $\Rc$,
and we state our main results about the special structures
of the block-Krylov subspaces associated with these two classes.
Proofs of these results are given in Section~\ref{sec-proofs} below.

\subsection{Case I.}

In this subsection, we assume that the matrices~\mr{KA0} are of the
form
\beq{MBA}{
\bea{rl}
\Mc &\!\!\!\! = \bigl(c \otimes I_{n_0}\bigr)
        \Matrix{M^{(1)} & M^{(2)} & \cdots & M^{(l)}}
    + \Sigma \otimes I_{n_0},\\[8pt]
\Rc &\!\!\!\! = c \otimes R,
\ea
}
where
\beq{MBB}{
\bea{rl}
M^{(i)} &\!\!\!\! \in \complex^{n_0 \times n_0},\q i=1,2,\ldots,l,\q
R \in \complex^{n_0\times m},\\[8pt]
c = \Matrix{c_1 \cr c_2 \cr \vdots \cr c_l} &\!\!\!\! \in \complex^l,\q
\mbox{and}\q 
\Sigma =\Matrix{\sigma_{11} & \sigma_{12} & \cdots & \sigma_{1l} \cr
                \sigma_{21} & \sigma_{22} & \cdots & \sigma_{2l} \cr
                \vdots & \vdots & & \vdots \cr
                \sigma_{l1} & \sigma_{l2} & \cdots & \sigma_{ll}}
    \in \complex^{l\times l}.
\ea
}
We assume that
\beq{MBC}{
c_i \not=0,\q i=1,2,\ldots,l.
}
Note that $\Mc \in \complex^{N \times N}$ and $\Rc \in \complex^{N \times m}$,
where
\beq{MBE}{
N := l n_0.
}

Our main result about the structure of the block-Krylov subspaces associated 
with the class of matrices~\mr{MBA} is as follows.
  
\begin{theorem} \label{theorem1}
Let\/~$\Mc$ and\/~$\Rc$ be matrices of the form\/~\mr{MBA} and\/~\mr{MBB},
and assume that\/~\mr{MBC} is satisfied.
Let\/~$\Vc \in \complex^{N\times N_0}$ be any basis matrix
of the block-Krylov subspaces induced by\/~$\Mc$ and\/~$\Rc$.
Then,\/~$\Vc$ can be represented in the form
\beq{MBG}{
\Vc = \Matrix{W \, U^{(1)}\cr \noalign{\vskip4pt}
              W \, U^{(2)}\cr \noalign{\vskip2pt}
                \vdots \cr \noalign{\vskip4pt}
              W \, U^{(l)}},
}
where $W\in \complex^{n_0 \times N_0}$ and, for each $i=1,2,\ldots,l$,
$U^{(i)} \in \complex^{N_0 \times N_0}$ is nonsingular and upper triangular.
\end{theorem}
 
The result of Theorem~\ref{theorem1} can be interpreted as follows.
Let
\beq{MXX}{
S_n := \Spa \Matrix{w_1 & w_2 & \cdots & w_n} \subseteq \complex^{n_0},\q 
n=1,2,\ldots,N_0,
}
denote the sequence of subspaces spanned by the leading columns of the
matrix
\beq{MXY}{
W = \Matrix{w_1 & w_2 & \cdots & w_{N_0}} \in \complex^{n_0 \times N_0}.
}
In view of~\mr{MBG}, for each $n=1,2,\ldots,N_0$, the $n$-th block-Krylov subspace
$\Kc_n(\Mc,\Rc)$, even though it is a subspace in $\complex^N$, consists of
$l$ `copies' of the same subspace $S_n$, which is a subspace of only $\complex^{n_0}$, 
where, by~\mr{MBE}, $n_0 = N/l$.
We stress that, in general, $S_n$ is not a block-Krylov subspace.

\subsection{Case II.}

In this subsection, we assume that the matrices~\mr{KA0} are of the 
form
\beq{MAA}{
\bea{rl}
\Mc &\!\!\!\! = \Matrix{C^{(1)} \cr \noalign{\vskip4pt}
                        C^{(2)} \cr \noalign{\vskip2pt}
                         \vdots \cr \noalign{\vskip4pt}
                       C^{(l)}}
        \Matrix{M^{(1)} & M^{(2)} & \cdots & M^{(l)}}
    + \Matrix{\sigma_1 I_{n_1}  & 0 & \cdots & 0 \cr
               0 & \sigma_2 I_{n_2} & \ddots & \vdots \cr
             \vdots & \ddots  & \ddots  & 0 \cr
                  0 & \cdots & 0 & \sigma_l I_{n_l}},\\[35pt]
\Rc &\!\!\!\! =  \Matrix{C^{(1)} \cr \noalign{\vskip4pt}
                        C^{(2)} \cr \noalign{\vskip2pt}
                         \vdots \cr \noalign{\vskip4pt}
                       C^{(l)}} R, 
\ea
}
where
\beq{MAB}{
\bea{rl}
& C^{(i)}  \in \complex^{n_i \times n_0},\q
M^{(i)} \in \complex^{n_0 \times n_i},\q \sigma_i \in \complex,\q 
            i=1,2,\ldots,l,\\[8pt]
& \mbox{and}\q R \in \complex^{n_0\times m}.
\ea
}
Note that $\Mc \in \complex^{N \times N}$ and $\Rc \in \complex^{N \times m}$,
where
\[
N := n_1 + n_2 + \cdots + n_l.
\]

Our main result about the structure of the block-Krylov subspaces associated 
with the class of matrices~\mr{MAA} is as follows.

\begin{theorem} \label{theorem2}
Let\/~$\Mc$ and\/~$\Rc$ be matrices of the form\/~\mr{MAA} and\/~\mr{MAB}.
Let\/~$\Vc \in \complex^{N\times N_0}$ be any basis matrix
of the block-Krylov subspaces induced by\/~$\Mc$ and\/~$\Rc$.
Then,\/~$\Vc$ can be represented in the form
\beq{MAG}{
\Vc = \Matrix{C^{(1)}\, W\, U^{(1)}\cr \noalign{\vskip4pt}
                C^{(2)}\, W\, U^{(2)}\cr \noalign{\vskip2pt}
                \vdots \cr \noalign{\vskip4pt}
                C^{(l)}\, W\, U^{(l)}},
}
where $W\in \complex^{n_0 \times N_0}$ and, for each $i=1,2,\ldots,l$,
$U^{(i)} \in \complex^{N_0 \times N_0}$ is nonsingular and upper triangular.
\end{theorem}

The result of Theorem~\ref{theorem2} can be interpreted as follows.
Let $S_n \subseteq \complex^{n_0}$, $n=1,2,\ldots,N_0$, again denote the sequence 
of subspaces spanned by the leading columns of the matrix $W$; as defined in~\mr{MXX}
and~\mr{MXY}.
In view of~\mr{MAG}, for each $n=1,2,\ldots,N_0$, the $n$-th block-Krylov subspace
$\Kc_n(\Mc,\Rc)$, even though it is a subspace in $\complex^N$, consists of
$l$ `copies' of the $C^{(i)}$-multiples, $i=1,2,\ldots,l$, of
the same subspace $S_n$, which is a subspace of $\complex^{n_0}$.

\section{Proofs} \label{sec-proofs}

In this section, we present proofs of Theorems~\ref{theorem1} 
and~\ref{theorem2}.

\subsection{Proof of Theorem~\ref{theorem1}.} \label{ssec-thm1}

Let $\Vc$ be a given basis matrix of the block Krylov subspaces induced
by $\Mc$ and $\Rc$.
We need to show that there exists a matrix $W$ and
nonsingular upper triangular matrices $U^{(i)}$, $i=1,2,\ldots,l$, such 
that~\mr{MBG} holds true.

Recall that any two basis matrices are connected by a relation of
the form~\mr{KAM}, where $\Uc$ is a nonsingular and upper triangular matrix.
Therefore, without loss of generality, we may assume that 
\beq{VVV}{
\Vc = \Vc(\Mc,\Rc)
}
is the particular basis matrix defined in~\mr{KAC}.
Furthermore, we partition any possible candidate matrices $W$ and $U^{(i)}$, 
$i=1,2,\ldots,l$, according to the block sizes of $\Vc(\Mc,\Rc)$
in~\mr{KAC}.
More precisely, we set
\beq{PBC}{
\bea{rl}
W &\!\!\!\! = \Matrix{W_1 & W_2 & \cdots & W_{k_0}},\\[16pt]
U^{(i)} &\!\!\!\! = \Matrix{U_{11}^{(i)} & U_{12}^{(i)} & \cdots & U_{1k_0}^{(i)}
                      \cr  \noalign{\vskip2pt}
                 0 & U_{22}^{(i)} & \ddots & \vdots
                      \cr \noalign{\vskip2pt}
               \vdots & \ddots & \ddots & \vdots 
                      \cr \noalign{\vskip4pt}
                 0 & \cdots & 0 & U_{k_0k_0}^{(i)}}, \q i=1,2,\ldots,l,
\ea
}
with subblocks $W_k \in \complex^{n_0 \times m_k}$ and nonsingular upper triangular
diagonal blocks $U_{kk}^{(i)} \in \complex^{m_k \times m_k}$ for all $i=1,2,\ldots,l$
and $k=1,2,\ldots,k_0$.
Inserting~\mr{VVV} and~\mr{PBC} into~\mr{MBG}, it follows that the
desired relation~\mr{MBG} holds true if, and only if, 
\beq{PBH}{
\Mc^{k-1} \Rc_k = \Matrix{\sum_{j=1}^k W_j U_{jk}^{(1)}
                                        \cr \noalign{\vskip4pt}
                          \sum_{j=1}^k W_j U_{jk}^{(2)}
                                        \cr \noalign{\vskip2pt}
                            \vdots \cr \noalign{\vskip4pt}
                          \sum_{j=1}^k W_j U_{jk}^{(l)}},\q
k=1,2,\ldots,k_0.
}
Therefore, it remains to construct the subblocks in~\mr{PBC} such that~\mr{PBH}
is satisfied.  To this end, we define these subblocks recursively as follows.

For each $k=1,2,\ldots,k_0$, we set
\beq{PBT}{ 
U_{jk}^{(i)} := \cases{
                  \displaystyle{\left(\sum_{t=1}^l 
                              \sigma_{i,t} U_{j,k-1}^{(t)}\right) E_k},  
                         & \mbox{for $j=1,2,\ldots,k-1$},\cr\noalign{\vskip12pt}
                   c_i I_{m_k}, & \mbox{for $j=k$},
                 }
}
for all $i=1,2,\ldots,l$, and
\beq{PBU}{
W_k := \cases{
            R E_1, & \mbox{if $k=1$},\cr\noalign{\vskip12pt}
          \displaystyle{\left(\sum_{i=1}^l M^{(i)} 
      \Biggl(\sum_{j=1}^{k-1} W_j U_{j,k-1}^{(i)} \Biggr) \right) E_k}, 
              & \mbox{if $k>1$}.
                 }
}
Here, the matrices $E_k$ are the ones from~\mr{KAD}.
We remark that, in view of assumption~\mr{MBC}, the subblocks $U^{(i)}_{kk}$ 
in~\mr{PBT} are all nonsingular.
Moreover, they are all diagonal and thus, in particular, upper triangular.

Using induction on $k$, we now show that the subblocks~\mr{PBT} and~\mr{PBU}
indeed satisfy~\mr{PBH}.
Recall from~\mr{MBA} that $\Rc = c\otimes R$ and from~\mr{KAD} (for $k=1$) 
that $\Rc_1 = \Rc E_1$. 
Together with the definitions of $U^{(i)}_{11}$, $i=1,2,\ldots,l$,
in~\mr{PBT} and of $W_1$ in~\mr{PBU}, it follows that
\[
\Rc_1 = c \otimes \bigl(R E_1\bigr) = \Matrix{c_1 R E_1 \cr \noalign{\vskip2pt}
                                  c_2 R E_1 \cr \noalign{\vskip2pt}
                                     \vdots \cr%
                                  c_l R E_1}
= \Matrix{W_1 U^{(1)}_{11} \cr \noalign{\vskip2pt}
          W_1 U^{(2)}_{11} \cr%
             \vdots \cr \noalign{\vskip2pt}
          W_1 U^{(l)}_{11}}.
\]
This is just~\mr{PBH} for $k=1$.
Let $1 < k \leq k_0$ and assume that~\mr{PBH} holds true for $k-1$.
Then, by multiplying the relation~\mr{PBH} (with $k$ replaced by $k-1$)
from the left by the matrix $\Mc$ from~\mr{MBA}, it follows that
\[
\bea{rl}
\Mc^{k-1} \Rc_{k-1} &\!\!\!\! = \Mc \Bigl(\Mc^{k-2} \Rc_{k-1} \Bigr) \\[8pt]
= \bigl(c \otimes I_{n_0}\bigr) &\!\!\!\! \displaystyle{\left(\sum_{i=1}^l M^{(i)} 
     \Biggl(\sum_{j=1}^{k-1} W_j U_{j,k-1}^{(i)} \Biggr) \right)} 
 + \Matrix{\sum_{t=1}^l \sigma_{1,t} \sum_{j=1}^{k-1} W_j U_{j,k-1}^{(t)}
                                                        \cr \noalign{\vskip4pt}
          \sum_{t=1}^l \sigma_{2,t} \sum_{j=1}^{k-1} W_j U_{j,k-1}^{(t)}
                                                        \cr \noalign{\vskip2pt}
          \vdots \cr \noalign{\vskip4pt}
          \sum_{t=1}^l \sigma_{l,t} \sum_{j=1}^{k-1} W_j U_{j,k-1}^{(t)}}.
\ea
\]
Multiplying this relation from the right by the
matrix $E_k$ from~\mr{KAD} and using the definitions of $U^{(i)}_{jk}$, 
$i=1,2,\ldots,l$, $j=1,2,\ldots,k$, in~\mr{PBT} and of $W_k$ in~\mr{PBU},
we obtain
\[
\bea{rl}
\Mc^{k-1} \Rc_k &\!\!\!\! = \Mc^{k-1} \Rc_{k-1} E_k\\[8pt]
&\!\!\!\! = \Matrix{W_k U_{kk}^{(1)} \cr \noalign{\vskip4pt}
                     W_k U_{kk}^{(2)} \cr \noalign{\vskip2pt}
          \vdots \cr \noalign{\vskip4pt}
                     W_k U_{kk}^{(l)}}
 + \Matrix{\sum_{j=1}^{k-1} W_j U_{jk}^{(1)}
                                        \cr \noalign{\vskip4pt}
                          \sum_{j=1}^{k-1} W_j U_{jk}^{(2)}
                                        \cr \noalign{\vskip2pt}
                            \vdots \cr \noalign{\vskip4pt}
                          \sum_{j=1}^{k-1} W_j U_{jk}^{(l)}}
           = \Matrix{\sum_{j=1}^k W_j U_{jk}^{(1)}
                                        \cr \noalign{\vskip4pt}
                          \sum_{j=1}^k W_j U_{jk}^{(2)}
                                        \cr \noalign{\vskip2pt}
                            \vdots \cr \noalign{\vskip4pt}
                          \sum_{j=1}^k W_j U_{jk}^{(l)}}.
\ea 
\]
This is just the desired relation~\mr{PBH}, and 
thus the proof of Theorem~\ref{theorem1} is complete.

\subsection{Proof of Theorem~\ref{theorem2}.} \label{ssec-thm2}

We proceed in the same fashion as in Subsection~\ref{ssec-thm1}.
Again, without loss of generality, we assume that the basis matrix $\Vc$
in~\mr{MAG} is given by~\mr{VVV}, and we partition the matrices $W$ and $U^{(i)}$,
$i=1,2,\ldots,l$, as in~\mr{PBC}.
Inserting~\mr{VVV} and~\mr{PBC} into~\mr{MAG}, it follows that the
desired relation~\mr{MAG} holds true if, and only if,
\beq{PAH}{
\Mc^{k-1} \Rc_k = \Matrix{C^{(1)}\, \sum_{j=1}^k W_j U_{jk}^{(1)}
                                        \cr \noalign{\vskip4pt}
                C^{(2)}\, \sum_{j=1}^k W_j U_{jk}^{(2)}
                                        \cr \noalign{\vskip2pt}
                \vdots \cr \noalign{\vskip4pt}
                C^{(l)}\, \sum_{j=1}^k W_j U_{jk}^{(l)}},\q
k=1,2,\ldots,k_0.
}
Therefore, it remains to construct the subblocks in~\mr{PBC} such that~\mr{PAH}
is satisfied.  To this end, we define these subblocks recursively as follows.
 
For $k=1,2,\ldots,k_0$, we set
\beq{PAT}{ 
U_{jk}^{(i)} := \cases{
                   \sigma_i U_{j,k-1}^{(i)} E_k,  
                         & \mbox{for $j=1,2,\ldots,k-1$},\cr\noalign{\vskip8pt}
                   I_{m_k}, & \mbox{for $j=k$},
                 }
}
for all $i=1,2,\ldots,l$, and 
\beq{PAU}{
W_k := \cases{
            R E_1, & \mbox{if $k=1$},\cr\noalign{\vskip12pt}
          \displaystyle{\left(\sum_{i=1}^l M^{(i)} C^{(i)} 
      \Biggl(\sum_{j=1}^{k-1} W_j U_{j,k-1}^{(i)} \Biggr) \right) E_k}, 
              & \mbox{if $k>1$}.
                 }
}
Here, again, the matrices $E_k$ are the ones from~\mr{KAD}.
 
Using induction on $k$, we now show that the subblocks~\mr{PAT} and~\mr{PAU}
indeed satisfy~\mr{PAH}.
Recall that $\Rc$ is of form~\mr{MAA} and that, by~\mr{KAD} (for $k=1$),
$\Rc_1 = \Rc E_1$. 
Together with the definitions of $U^{(i)}_{11}$, $i=1,2,\ldots,l$,
in~\mr{PAT} and of $W_1$ in~\mr{PAU}, it follows that
\[
\Rc_1 = \Matrix{C^{(1)} \cr \noalign{\vskip4pt}
                        C^{(2)} \cr \noalign{\vskip2pt}
                         \vdots \cr \noalign{\vskip4pt}
                       C^{(l)}} \bigl(R E_1\bigr) 
 = \Matrix{C^{(1)}\, W_1\, U_{jk}^{(1)}
                      \cr \noalign{\vskip4pt}
                C^{(2)}\, W_1\, U_{11}^{(2)}
                                        \cr \noalign{\vskip2pt}
                \vdots \cr \noalign{\vskip4pt}
                C^{(l)}\, W_1\, U_{11}^{(l)}}.
\]
This is just~\mr{PAH} for $k=1$.
Let $1 < k \leq k_0$ and assume that~\mr{PAH} holds true for $k-1$.
Then, by multiplying the relation~\mr{PAH} (with $k$ replaced by $k-1$)
from the left by the matrix $\Mc$ from~\mr{MAA}, it follows that
\[
\bea{rl}
\Mc^{k-1} \Rc_{k-1} &\!\!\!\! = \Mc \Bigl(\Mc^{k-2} \Rc_{k-1} \Bigr) \\[8pt]
=  \Matrix{C^{(1)} \cr \noalign{\vskip4pt}
                        C^{(2)} \cr \noalign{\vskip2pt}
                         \vdots \cr \noalign{\vskip4pt}
                       C^{(l)}}
        &\!\!\!\! \displaystyle{\left(\sum_{i=1}^l M^{(i)} 
    C^{(i)} \Biggl(\sum_{j=1}^{k-1} W_j U_{j,k-1}^{(i)} \Biggr) \right)} 
 + \Matrix{C^{(1)} \sum_{j=1}^{k-1} W_j \sigma_1 U_{j,k-1}^{(1)}
                                                        \cr \noalign{\vskip4pt}
           C^{(2)} \sum_{j=1}^{k-1} W_j \sigma_2 U_{j,k-1}^{(2)}
                                                        \cr \noalign{\vskip2pt}
          \vdots \cr \noalign{\vskip4pt}
           C^{(l)} \sum_{j=1}^{k-1} W_j \sigma_l U_{j,k-1}^{(l)}}.
\ea
\]
Multiplying this relation from the right by the
matrix $E_k$ from~\mr{KAD} and using the definitions of $U^{(i)}_{jk}$, 
$i=1,2,\ldots,l$, $j=1,2,\ldots,k$, in~\mr{PAT} and of $W_k$ in~\mr{PAU},
we obtain
\[ 
\bea{rl}
\Mc^{k-1} \Rc_k &\!\!\!\! = \Mc^{k-1} \Rc_{k-1} E_k\\[8pt]
 = &\!\!\!\! \Matrix{C^{(1)} W_k U_{kk}^{(1)} \cr \noalign{\vskip4pt}
                        C^{(2)} W_k U_{kk}^{(2)}\cr \noalign{\vskip2pt}
                         \vdots \cr \noalign{\vskip4pt}
                       C^{(l)} W_k U_{kk}^{(l)}}
 + \Matrix{C^{(1)} \sum_{j=1}^{k-1} W_j  U_{jk}^{(1)}
                                                        \cr \noalign{\vskip4pt}
           C^{(2)} \sum_{j=1}^{k-1} W_j  U_{jk}^{(2)}
                                                        \cr \noalign{\vskip2pt}
          \vdots \cr \noalign{\vskip4pt}
           C^{(l)} \sum_{j=1}^{k-1} W_j  U_{jk}^{(l)}}
=  \Matrix{C^{(1)} \sum_{j=1}^{k} W_j U_{jk}^{(1)}
                                                        \cr \noalign{\vskip4pt}
           C^{(2)} \sum_{j=1}^{k} W_j U_{jk}^{(2)}
                                                        \cr \noalign{\vskip2pt}
          \vdots \cr \noalign{\vskip4pt}
           C^{(l)} \sum_{j=1}^{k} W_j U_{jk}^{(l)}}.
\ea 
\]
This is just the desired relation~\mr{PAH}, and 
thus the proof of Theorem~\ref{theorem2} is complete.

\section{Matrices arising in higher-order linear dynamical systems.} 
\label{sec-higher}

In this section, we show how block-Krylov subspaces $\Kc_n(\Mc,\Rc)$ with 
matrices $\Mc$ and $\Rc$ of the form~\mr{MBA} arise in the
context of higher-order linear dynamical systems.

\subsection{General time-invariant linear dynamical systems.}

We consider general higher-order multi-input multi-output time-invariant 
linear dynamical systems.
We denote by $m$ and $p$ the number of inputs and outputs,
respectively, and by $l$ the order of such systems.
In the following, the only assumption on $m$, $p$, and $l$ is
that $m$, $p$, $l \geq 1$.

An \textit{$m$-input $p$-output time-invariant linear dynamical 
system of order $l$} is a system of differential-algebraic 
equations (DAEs) of the following form:
\beq{AAA}{
\bea{rl} 
& \quad\ P_l \displaystyle{\frac{d^l}{dt^l}} x(t)
 + P_{l-1} \displaystyle{\frac{d^{l-1}}{dt^{l-1}}} x(t)
 + \cdots + P_1 \displaystyle{\frac{d}{dt}} x(t) + P_0 x(t)
                   = B u(t), \\[12pt]
& y(t) = D u(t) 
  + L_{l-1} \displaystyle{\frac{d^{l-1}}{dt^{l-1}}} x(t)
 + \cdots + L_1 \displaystyle{\frac{d}{dt}} x(t) + L_0 x(t).
\ea 
} 
Here, $P_i\in \complex^{n_0 \times n_0}$, $0\leq i \leq l$,
$B \in \complex^{n_0\times m}$, $D \in \complex^{p\times m}$,
and $L_j \in \complex^{p\times n_0}$, $0\leq j < l$, 
are given matrices, and $n_0$ is called the state-space dimension
of~\mr{AAA}.
Moreover, in~\mr{AAA},
 $u: [t_0, \infty) \mapsto \complex^m$
is a given input function, $t_0 \in \real$ is a given
initial time, the components of the vector-valued
function $x: [t_0, \infty) \mapsto \complex^{n_0}$ are the so-called
state variables, and $y: [t_0, \infty) \mapsto \complex^p$
is the output function.
The system is completed by initial conditions of the form
\beq{AAB}{
\displaystyle{\frac{d^j}{dt^j}} x(t) \biggm|_{t=t_0} 
= x_0^{(j)},\q 0\leq j < l,
}
where $x_0^{(j)} \in \complex^{n_0}$, $0\leq j < l$, are given
vectors.

We stress that the matrix $P_l$ is allowed to be singular,
and thus the first equation in~\mr{AAA} is indeed a system
of DAEs in general.
Our only assumption on the matrices $P_i$, $0 \leq i \leq l$,
in~\mr{AAA} is that the $n_0\times n_0$-matrix-valued polynomial
\beq{AAE}{
P(s) := s^l P_l + s^{l-1} P_{l-1} + \cdots + s P_1 + P_0,\q s \in \complex,
}
is \textit{regular}, i.e.,
the matrix $P(s)$ is singular only for finitely many values 
of $s\in \complex$; see, e.g., \cite[Part II]{GohLR82}.

\subsection{Equivalent first-order formulation.} 

It is well known (see, e.g., \cite[Chapter 7]{GohLR82}) that 
any $l$-th-order system ~\mr{AAA} (with state-space dimension $n_0$)
is equivalent to a first-order system with state-space dimension $N := l n_0$.
Indeed, it is easy to verify that the $l$-th-order system~\mr{AAA}
with initial conditions~\mr{AAB} is equivalent to the
first-order system
\beq{CAA}{
\bea{rl} 
\Ec \displaystyle{\frac{d}{dt}} z(t) - \Ac z(t) 
      &\!\!\!\! = \Bc u(t), \\[12pt]
y(t)  &\!\!\!\! = \Dc u(t) + \Lc z(t), \\[12pt]
z(t_0) &\!\!\!\! = z_0,
\ea  
}
where
\beq{CAC}{
\bea{rl}
z(t) &\!\!\!\! := \Matrix{x(t)\cr 
                         \noalign{\vskip3pt}
                         \frac{d}{dt} x(t) \cr 
                         \vdots \cr 
                         \noalign{\vskip3pt}
                         \frac{d^{l-1}}{dt^{l-1}} x(t)},\q
z_0 := \Matrix{x_0^{(0)} \cr 
              \noalign{\vskip3pt}
	      x_0^{(1)} \cr
              \vdots \cr
	      \noalign{\vskip3pt}
              x_0^{(l-1)}},\\[42pt]
& \Lc := \Matrix{L_0 & L_1 & \cdots & L_{l-1}},\q
\Bc := \Matrix{0 \cr \vdots \cr 0 \cr \noalign{\vskip3pt} B},\q
\Dc := D,\\[38pt] 
\Ec  &\!\!\!\! := \Matrix{I  & 0 & 0 & \cdots & 0 \cr
               \noalign{\vskip2pt}
               0 & I & 0 & \cdots & 0\cr
              \vdots & \ddots & \ddots & \ddots & \vdots \cr
                  0 & \cdots & 0  & I  & 0 \cr
               \noalign{\vskip2pt}
                  0 & \cdots & 0 & 0 & P_l},\q	
\Ac := -\Matrix{   0 & \!\! -I & 0 & \cdots & 0\cr 
                  0 & 0 & \!\! -I & \ddots & \vdots \cr
             \vdots & \ddots &  \ddots  & \ddots  & 0 \cr
                  0 & \cdots & 0 & 0 & \!\! -I \cr
             \noalign{\vskip2pt}
 P_0 & P_1 & P_2 & \cdots & P_{l-1}},
\ea
} 
and $I=I_{n_0}$ is the $n_0\times n_0$ identity matrix.

It is easy to see that, for any given $s\in \complex$, the 
matrix $s\, \Ec - \Ac$ 
is singular if, and only if, the matrix $P(s)$ defined in~\mr{AAE}
is singular.
Therefore, our assumption on the regularity of the matrix polynomial~\mr{AAE}
is equivalent to the regularity of the matrix pencil $s\, \Ec - \Ac$.
This guarantees that the matrix $s\, \Ec - \Ac $
is singular only for finitely many values of $s\in \complex$, and
that
\beq{CAE}{
H(s) := \Dc + \Lc \bigl(s\, \Ec - \Ac\bigr)^{-1} \Bc,\q s \in \complex,
}
is a well-defined $p\times m$-matrix-valued rational function.
We remark that~\mr{CAE} is called the frequency-domain \textit{transfer function} 
of~\mr{CAA}.

\subsection{Pad\'e-type model reduction.} \label{ssec-pade1} 

A reduced-order model of~\mr{CAA} is a linear dynamical system
of the same type as~\mr{CAA}, but with reduced state-space dimension, say $n$, 
instead of the original state-space dimension $N$.
More precisely, a \textit{reduced-order model} of~\mr{CAA} with state-space
dimension $n$ is a system of the form
\beq{PAA}{
\bea{rl} 
\Ec_n \displaystyle{\frac{d}{dt}} \tilde{z}(t) - \Ac_n \tilde{z}(t) 
\      &\!\!\!\! = \Bc_n u(t), \\[12pt]
\tilde{y}(t)  &\!\!\!\! = \Dc_n u(t) + \Lc_n \tilde{z}(t), \\[12pt]
\tilde{z}(t_0) &\!\!\!\! = \tilde{z}_0,
\ea  
}
where $\Ac_n$, $\Ec_n \in \complex^{n\times n}$,
$\Bc_n \in \complex^{n\times m}$, $\Dc_n \in \complex^{p\times m}$,
$\Lc_n \in \complex^{p\times n}$, and $\tilde{z}_0 \in \complex^n$.
The problem of model reduction then is to construct data matrices 
$\Ac_n$, $\Ec_n$, $\Bc_n$, $\Dc_n$, and $\Lc_n$ such that~\mr{PAA}
is a good approximation of the original system~\mr{CAA}, even 
for $n\ll N$.

A possible approach, which is intimately related to block-Krylov subspaces,
is Pad\'e and Pad\'e-type model reduction; see, e.g.,~\cite{Fre03b,Fre05a} and
the references given there.
Let $s_0\in \complex$ be a suitably chosen expansion point, and in particular,
let $s_0$ be such that the matrix $s_0\, \Ec - \Ac$ is nonsingular. 
The reduced system~\mr{PAA} is said to be an $n$-th Pad\'e model of the original
system~\mr{CAA} if the reduced-order transfer function
\[
H_n(s) := \Dc_n + \Lc_n \bigl(s\, \Ec_n - \Ac_n\bigr)^{-1} \Bc_n,\q s \in \complex,
\]
and the original transfer function~\mr{CAE}, $H$, agree in as many leading
Taylor coefficients about the expansion point $s_0$ as possible, i.e.,
\beq{PAD}{
H_n(s)  = H(s) + \Oc\bigl((s-s_0)^{q(n)}\bigr),
}
where $q(n)$ is as large as possible.
While Pad\'e models are optimal in the sense of~\mr{PAD}, in general, they do 
not preserve other desirable properties of the original system.
Preserving such properties is often possible by relaxing~\mr{PAD} to
\beq{PTY}{
H_n(s)  = H(s) + \Oc\bigl((s-s_0)^{\hat{q}}\bigr),
} 
where $\hat{q} < q(n)$.
The reduced system~\mr{PAA} is said to be an $n$-th Pad\'e-type model of the original
system~\mr{CAA} if a property of the form~\mr{PTY} is satisfied.

Both $n$-th Pad\'e and Pad\'e-type models can be generated via Krylov-subspace
machinery; see, e.g.,~\cite{Fre03b,Fre05a} and the references given there.
To this end, the original transfer function~\mr{CAE} is rewritten in the form
\[
H(s) = \Dc + \Lc \bigl(\Ic + (s - s_0)\, \Mc\bigr)^{-1} \Rc,
\]
where  
\beq{EAC}{
\Mc := \bigl(s_0\, \Ec - \Ac\bigr)^{-1} \Ec\q \mbox{and}\q
\Rc := \bigl(s_0\, \Ec - \Ac\bigr)^{-1} \Bc.
}
Pad\'e-type models are then obtained by projecting the data matrices in~\mr{CAA}
onto the block-Krylov subspaces $\Kc_n(Mc,\Rc)$ induced by the matrices~\mr{EAC}.
Similarly, Pad\'e models can be generated via two-sided projections involving
the right and left block-Krylov subspaces $\Kc_n(\Mc,\Rc)$ and $\Kc_n(\Mc^H,\Lc^H)$.

\subsection{Structure of the matrices $\Mc$ and $\Rc$.} 

Recall that, in this section, we are concerned with general $l$-th-order systems of 
the form~\mr{AAA}.
In this case, the matrices $\Ac$, $\Ec$, and $\Bc$ in~\mr{EAC} are the
ones defined in~\mr{CAC}.
Furthermore, the expansion point $s_0\in \complex$ in~\mr{EAC} is such that
the matrix $s_0\, \Ec - \Ac$ is nonsingular, or, equivalently, the matrix
\beq{PSO}{
P(s_0) = s_0^l P_l + s_0^{l-1} P_{l-1} + \cdots + s_0 P_1 + P_0\q
\mbox{is nonsingular}.
}

Next, we set
\beq{MDEF}{
M^{(i)} := \bigl( P(s_0) \bigr)^{-1} \sum_{j=0}^{l-i} s_0^j P_{i+j},\q
i=1,2,\ldots,l,
} 
and
\beq{RDEF}{
R := \bigl(P(s_0) \bigr)^{-1} B.
}
Using the definitions of $\Ac$, $\Ec$, and $\Bc$ in~\mr{EAC}, together 
with~\mr{MDEF} and~\mr{RDEF}, one can show that the matrices~\mr{EAC}
have the representations
\beq{PIA}{
\bea{rl}
\Mc = &\!\!\!\! 
\Matrix{M^{(1)} & M^{(2)} & M^{(3)} & \cdots &  M^{(l)} \cr \noalign{\vskip2pt}
s_0 M^{(1)} & s_0 M^{(2)} & s_0 M^{(3)} & \cdots
                                 & s_0 M^{(l)} \cr \noalign{\vskip4pt}
s_0^2 M^{(1)} & s_0^2 M^{(2)} & s_0^2 M^{(3)} & \cdots &  s_0^2 M^{(l)} 
                                                \cr \noalign{\vskip4pt}
\vdots & \vdots &  \vdots &   & \vdots \cr \noalign{\vskip4pt}
s_0^{l-2} M^{(1)} & s_0^{l-2} M^{(2)} & s_0^{l-2} M^{(3)} & \cdots &  
                                  s_0^{l-2} M^{(l)}}\\[48pt]
   &\qquad \qquad - \Matrix{0 & 0 & \cdots & \cdots & 0 \cr
                    I_{n_0} & 0 & \ddots & & \vdots \cr
                    s_0 I_{n_0} & I_{n_0} & 0 & \ddots & \vdots \cr
                   \vdots & \ddots & \ddots & \ddots & \vdots \cr
                  s_0^{l-2} I_{n_0} & \cdots & s_0 I_{n_0} & I_{n_0} & 0},
\ea
}
and
\beq{PIC}{
\Rc = 
      \Matrix{I_{n_0} \cr
              s_0 I_{n_0} \cr \noalign{\vskip4pt}
              s_0^2 I_{n_0} \cr \noalign{\vskip4pt}
              \vdots \cr \noalign{\vskip4pt}
             s_0^{l-1} I_{n_0}} R.
}
Proofs of~\mr{PIA} and~\mr{PIC} are given in Appendix~A.

Note that the matrices $\Mc$ and $\Rc$ in~\mr{PIA} and~\mr{PIC} are 
a special instance of the class of matrices~\mr{MBA}, with $c$ and $\Sigma$
given by
\[
c:= \Matrix{1 \cr 
            s_0 \cr \noalign{\vskip2pt}
            s_0^2 \cr
            \vdots \cr
            s_0^{l-1}}
\q \mbox{and}\q
\Sigma := - \Matrix{0 & 0 & \cdots & \cdots & 0 \cr
                    1 & 0 & \ddots & & \vdots \cr
                    s_0 & 1 & 0 & \ddots & \vdots \cr
                   \vdots & \ddots & \ddots & \ddots & \vdots \cr
                    s_0^{l-2} & \cdots & s_0 & 1 & 0}.
\]
Furthermore, provided that $s_0\not=0$, the assumption on $c$ in~\mr{MBC}
is satisfied.
We remark that for the case $s_0=0$, $\Mc$ reduces
to a block companion matrix, and $\Rc$ reduces to a multiple of the first block 
unit vector.
We do not consider this case, which is fundamentally different from the 
case $s_0\not=0$, in this paper.
    
\section{Matrices arising in first-order integro-DAEs.}
\label{sec-second}
   
An important special case of~\mr{AAA} is 
\textit{second-order systems}, that is, $l=2$ in~\mr{AAA}.
For example, second-order systems arise in structural 
dynamics~\cite{Prz85,CraH88,SuC91}, circuit 
analysis~\cite[Chapter 3]{Tut77},
and computational electromagnetics~\cite{WitSW03}.
However, in some of these applications, a more suitable formulation
of such systems is as systems of first-order 
integro-differential-algebraic equations (integro-DAEs).
For example, this is the case for passive systems such as 
RCL electrical circuits consisting of only resistors,
capacitors, and inductors; see, e.g., 
\cite[Chapter 1]{LozBEM00}, \cite[Chapter 2]{DeC89}, and \cite{Fre04c,Fre05a}.
In this section, we show how block-Krylov subspaces $\Kc_n(\Mc,\Rc)$ with
matrices $\Mc$ and $\Rc$ of the form~\mr{MAA} arise in the
context of such systems of first-order 
integro-DAEs.
 
\subsection{Systems of first-order integro-DAEs.}

We consider $m$-input $p$-output systems of first-order 
\textit{integro-DAEs} of the following form: 
\beq{ASA}{ 
\bea{rl} 
P_1 \displaystyle{\frac{d}{dt}} x(t)
    + P_0 x(t) + P_{-1}
 \displaystyle{\int_{t_0}^t x(\tau)\, d\tau} &\!\!\!\! =
       B u(t), \\[12pt]
y(t)  &\!\!\!\! = D u(t) + L x(t), \\[12pt]
x(t_0) &\!\!\!\! = x_0. 
\ea 
}
Here, $P_{-1}$, $P_0$, $P_1 \in \complex^{n_0 \times n_0}$,
$B \in \complex^{n_0\times m}$, $D \in \complex^{p\times m}$,
and $L \in \complex^{p\times n_0}$ are given matrices,
$t_0 \in \real$ is a given initial time, 
and $x_0 \in \complex^{n_0}$ is a given vector of initial values.

We stress that the matrix $P_1$ is allowed to be singular,
and thus the first equation in~\mr{ASA} is indeed a system
of integro-DAEs in general.
Our only assumption on the matrices $P_{-1}$, $P_0$, and $P_1$
in~\mr{ASA} is that the $n_0 \times n_0$-matrix-valued rational
function
\[
Q(s) := s P_1 + P_0 + \sinv P_{-1},\q s\in \complex,
\]
is \textit{regular}, i.e., the matrix $Q(s)$ is singular only 
for finitely many values of $s\in \complex$.  

In practical applications, the matrices $P_0$ and $P_1$ are
usually sparse, while the matrix $P_{-1}$ is not always sparse.
However, in those cases where the matrix $P_{-1}$ itself is
dense, $P_{-1}$ is given as a product of the form
\beq{AF1}{
P_{-1} = F_1 G F_2^H
}
or
\beq{AF2}{
P_{-1} = F_1 G^{-1} F_2^H,\q \mbox{with  nonsingular $G$},
}
where $F_1$, $F_2 \in \complex^{n_0\times \hat{n}_0}$ and
$G\in \complex^{\hat{n}_0 \times \hat{n}_0}$ are
sparse matrices.
We stress that in the case~\mr{AF1}, the matrix $G$ is
not required to be nonsingular.
In particular, for any matrix $P_{-1} \in \complex^{n_0\times n_0}$,
there is always the trivial factorization~\mr{AF1} with
$F_1 = F_2 = I_{n_0}$ and $G = P_{-1}$.
Therefore, in the following, we assume that the matrix $P_{-1}$ 
in~\mr{ASA} is given by a product of the form~\mr{AF1}
or~\mr{AF2}.

\subsection{Equivalent first-order formulations.} 

In analogy to the case of higher-order systems~\mr{AAA},
any system of integro-DAEs of the form~\mr{ASA}
is equivalent to a first-order system of the form~\mr{CAA}.
In this subsection, we present such equivalent first-order formulations.

We distinguish the two cases~\mr{AF1} and~\mr{AF2}.
First assume that $P_{-1}$ is given by~\mr{AF1}.
In this case, we set 
\beq{CAG}{ 
   z_1(t) := x(t)\q \mbox{and}\q
   z_2(t) := F_2^H \displaystyle{\int_{t_0}^t x(\tau)\, d\tau}.
}
By~\mr{AF1} and~\mr{CAG}, the first relation in~\mr{ASA} can
be rewritten as follows:
\beq{CAI}{
P_1 z^{\prime}_1(t) + P_0 z_1(t) + F_1 G z_2(t) = B u(t).
}
Moreover,~\mr{CAG} implies that
\beq{CAK}{
z^{\prime}_2(t) = F_2^H z_1(t).
}
It follows from~\mr{CAG}--\mr{CAK} that the system of 
integro-DAEs~\mr{ASA} (with $P_{-1}$ given by~\mr{AF1}) is equivalent 
to a first-order system~\mr{CAA} where
\beq{CAM}{
\bea{rl}
z(t) &\!\!\!\! := \Matrix{z_1(t)\cr \noalign{\vskip2pt}
                        z_2(t)},\q
z_0 := \Matrix{x_0\cr \noalign{\vskip2pt} 0},\q
\Lc := \Matrix{L & 0},\q
\Bc := \Matrix{B \cr 0},\\[15pt]
\Dc &\!\!\!\! := D,\q
\Ac := \Matrix{-P_0 & -F_1 G \cr \noalign{\vskip2pt}
                F_2^H & 0},\q
\Ec := \Matrix{P_1 & 0 \cr \noalign{\vskip2pt}
                 0 & I_{\hat{n}_0}}.
\ea
} 

Next, we assume that $P_{-1}$ is given by~\mr{AF2}.
In this case, we set
\beq{JQA}{
   z_1(t) := x(t)\q \mbox{and}\q 
   z_2(t) := G^{-1} F_2^H \displaystyle{\int_{t_0}^t x(\tau)\, d\tau}.
}
By~\mr{AF2} and~\mr{JQA}, the first relation in~\mr{ASA} can
be rewritten as follows:
\beq{JQB}{
P_1 z^{\prime}_1(t) + P_0 z_1(t) + F_1 z_2(t) = B u(t).
}
Moreover,~\mr{JQA} implies that
\beq{JQC}{ 
G z^{\prime}_2(t) = F_2^H z_1(t). 
} 
It follows from~\mr{JQA}--\mr{JQC} that the system of 
integro-DAEs~\mr{ASA} (with $P_{-1}$ given by~\mr{AF2}) is 
equivalent to a first-order system~\mr{CAA} where
\beq{CAO}{
\bea{rl}
z(t) &\!\!\!\! := \Matrix{z_1(t)\cr \noalign{\vskip2pt}
                     z_2(t)},\q
z_0 := \Matrix{x_0\cr  \noalign{\vskip2pt} 0},\q
\Lc := \Matrix{L & 0},\q
\Bc := \Matrix{B \cr  \noalign{\vskip2pt} 0},\\[15pt]
\Dc &\!\!\!\! := D,\q
\Ac := \Matrix{-P_0 & -F_1 \cr \noalign{\vskip2pt}
               F_2^H & 0},\q
\Ec := \Matrix{P_1 & 0 \cr \noalign{\vskip2pt}
                 0 & G}.
\ea 
}

\subsection{Pad\'e and Pad\'e-type model reduction.} \label{ssec-pade2}

Just as in Subsection~\ref{ssec-pade1}, based on the equivalent first-order 
formulations defined in~\mr{CAM}, respectively~\mr{CAO}, one can again introduce the 
notion of Pad\'e and Pad\'e-type reduced-order models of systems of 
integro-DAEs~\mr{ASA}.
In this case, we assume that the expansion point $s_0\in \complex$ is chosen
such that $s_0\not = 0$ and the matrix 
\beq{QS0}{
Q_0 := Q(s_0) = s_0 P_1 + P_0 + \displaystyle{\szinv} P_{-1}
}
is nonsingular.
One readily verifies that this condition is equivalent to the nonsingularity of
the matrix $s_0\, \Ec - \Ac$.
The matrices that induce the relevant block-Krylov subspaces $\Kc_n(\Mc,\Rc)$
for Pad\'e and Pad\'e-type model reduction are again given by
\beq{EAC2}{
\Mc := \bigl(s_0\, \Ec - \Ac\bigr)^{-1} \Ec\q \mbox{and}\q
\Rc := \bigl(s_0\, \Ec - \Ac\bigr)^{-1} \Bc,
}
where $\Ac$, $\Ec$, and $\Bc$ are now the matrices defined in~\mr{CAM},
respectively~\mr{CAO}.

\subsection{Structure of the matrices $\Mc$ and $\Rc$.}

In this subsection, we describe the structure of the matrices $\Mc$ and $\Rc$.

Again, we distinguish the two cases~\mr{AF1} and~\mr{AF2}.
First assume that $P_{-1}$ is given by~\mr{AF1}.
Using the definitions of $\Ac$, $\Ec$, $\Bc$ in~\mr{CAM}, and of $Q_0$ in~\mr{QS0},
one can show that the matrices~\mr{EAC2} 
have the representations
\beq{STA}{
\bea{rl}
\Mc &\!\!\!\! = \Matrix{I_{n_0} \cr \noalign{\vskip4pt}
              \szinv F_2^H}
      \Matrix{Q_0^{-1} P_1 & -\szinv  Q_0^{-1} F_1 G}
      + \Matrix{0 & 0 \cr \noalign{\vskip4pt} 
                0 & \szinv I_{\hat{n}_0}},\\[16pt] 
\Rc &\!\!\!\! = \Matrix{I_{n_0} \cr \noalign{\vskip4pt}
              \szinv F_2^H} Q_0^{-1} B.
\ea
}
The matrices $\Mc$ and $\Rc$ in~\mr{STA} are 
a special instance of the class of matrices~\mr{MAA}, with the integers and
matrices~\mr{MAB} chosen as follows:
\[
\bea{rl}
l  &\!\!\!\! := 2, \q n_1 := n_0,\q n_2 := \hat{n}_0,\\[8pt]
C^{(1)} &\!\!\!\! := I_{n_0},\q 
   C^{(2)} := \displaystyle{\szinv} F_2^H,\\[8pt]
M^{(1)} &\!\!\!\! := Q_0^{-1} P_1,\q 
   M^{(2)} := - \displaystyle{\szinv} Q_0^{-1} F_1 G,\\[8pt]
\sigma_1 &\!\!\!\! := 0,\q \sigma_2 := \displaystyle{\szinv},\q R := Q_0^{-1} B.
\ea
\]

Next, we assume that $P_{-1}$ is given by~\mr{AF2}.
Using the definitions of $\Ac$, $\Ec$, $\Bc$ in~\mr{CAO}, and of $Q_0$ in~\mr{QS0},
one can show that the matrices~\mr{EAC2} 
have the representations
\beq{STM}{
\bea{rl}
\Mc &\!\!\!\! = \Matrix{I_{n_0} \cr \noalign{\vskip4pt}
              \szinv G^{-1} F_2^H}
      \Matrix{Q_0^{-1} P_1 & -\szinv  Q_0^{-1} F_1}
      + \Matrix{0 & 0 \cr \noalign{\vskip4pt} 
                0 & \szinv I_{\hat{n}_0}},\\[16pt]
\Rc &\!\!\!\! = \Matrix{I_{n_0} \cr \noalign{\vskip4pt}
              \szinv G^{-1} F_2^H} Q_0^{-1} B.
\ea
}
Note that the matrices $\Mc$ and $\Rc$ in~\mr{STM} are 
a special instance of the class of matrices~\mr{MAA}, with the integers and
matrices~\mr{MAB} chosen as follows:
\[
\bea{rl}
l  &\!\!\!\! := 2, \q n_1 := n_0,\q n_2 := \hat{n}_0,\\[8pt]
C^{(1)} &\!\!\!\! := I_{n_0},\q 
   C^{(2)} := \displaystyle{\szinv} G^{-1} F_2^H,\\[8pt]
M^{(1)} &\!\!\!\! := Q_0^{-1} P_1,\q 
   M^{(2)} := - \displaystyle{\szinv} Q_0^{-1} F_1,\\[8pt]
\sigma_1 &\!\!\!\! := 0,\q \sigma_2 := \displaystyle{\szinv},\q R := Q_0^{-1} B.
\ea
\]

Proofs of~\mr{STA} and~\mr{STM} are given in Appendix B.

\section{Concluding remarks.} \label{sec-concluding}
 
We have introduced two classes of structured matrices, which 
include the matrices of first-order
formulations of higher-order linear dynamical systems as special
cases.
As our main results, we have shown that the block-Krylov subspaces 
induced by the matrices in theses classes exhibit special structures.
Roughly speaking, for both classes, the associated structured 
block-Krylov subspaces consist of multiple copies of certain 
subspaces of the state space of the original higher-order system.
Note that the dimension of the state space of the first-order
formulation is $l$ times the dimension of the original $l$-th-order system.
Our results show that in order to construct basis vectors for the
block-Krylov subspaces of the higher-dimensional first-order state-space,
it is sufficient to construct basis vectors for certain subspaces of the 
lower-dimensional $l$-th-order state space.
The problem of the efficient and numerically stable construction 
of basis vectors of these subspaces is beyond the scope of this paper.
Such algorithms will be described in a forthcoming report.

\baselineskip=0.9\normalbaselineskip
  
{\small
\bibliographystyle{bit} 
\bibliography{../siam_book/biblio}
}

\section*{Appendix A.}

In this appendix, we establish the representations~\mr{PIA} and~\mr{PIC}.
To this end, we set
\beq{PDEF}{
\hat{P}_i:= \sum_{j=0}^{l-i} s_0^j P_{i+j},\q
i=0,1,\ldots,l.
}
In view of~\mr{PSO} and~\mr{MDEF}, we then have
\beq{XAA1}{
M^{(i)} := \hat{P}_0^{-1} \hat{P}_i,\q
i=1,2,\ldots,l, \q \mbox{and}\q 
\hat{P}_0 = P(s_0).
}
Using the definitions of $\Ac$ and $\Ec$ in~\mr{CAC}, as well as~\mr{PDEF},
one readily verifies that 
\beq{XAA2}{
\bea{rl}
s_0\, \Ec - \Ac &\!\!\!\! =
 \Matrix{s_0 I & -I & 0 & \cdots & 0\cr
                  0 & s_0 I & -I & \ddots & \vdots \cr
             \vdots & \ddots &  \ddots  & \ddots  & 0 \cr
                  0 & \cdots & 0 & s_0 I & -I\cr
             \noalign{\vskip2pt}
    P_0 & P_1 & \cdots & P_{l-2} & P_{l-1} + s_0 P_l}\\[40pt]
&\!\!\!\! =
       \Matrix{   0 & \!\! -I & 0 & \cdots & 0\cr 
                  0 & 0 & \!\! -I & \ddots & \vdots \cr 
             \vdots & \ddots &  \ddots  & \ddots  & 0 \cr
                  0 & \cdots & 0 & 0 & \!\! -I \cr
             \noalign{\vskip2pt}
    \hat{P}_0 & \hat{P}_1 & \hat{P}_2 & \cdots & 
                 \hat{P}_{l-1}} 
\Matrix{I & 0 & 0 & \cdots & 0 \cr
               -s_0 I & I & 0 & \cdots & 0\cr
                  0 & -s_0 I & I & \ddots & \vdots \cr
             \vdots & \ddots &  \ddots  & \ddots  & 0 \cr
                  0 & \cdots & 0 & -s_0 I & I}.
\ea
}
Note that, in view of~\mr{XAA1}, we have
\[
       \Matrix{   0 & \!\! -I & 0 & \cdots & 0\cr 
                  0 & 0 & \!\! -I & \ddots & \vdots \cr 
             \vdots & \ddots &  \ddots  & \ddots  & 0 \cr
                  0 & \cdots & 0 & 0 & \!\! -I \cr
             \noalign{\vskip2pt}
    \hat{P}_0 & \hat{P}_1 & \hat{P}_2 & \cdots & 
                 \hat{P}_{l-1}}^{-1}
=  \Matrix{M^{(1)} & M^{(2)} & \cdots & M^{(l-1)} & \hat{P}_0^{-1} \cr
          \noalign{\vskip2pt}
          -I & 0 & 0 & \cdots & 0 \cr
          0 & -I & 0 & \cdots & 0 \cr
     \vdots & \ddots & \ddots & \ddots & \vdots \cr
     0 & \cdots & 0 & -I  & 0}.
\]
By inverting the two factors on the right-hand side of~\mr{XAA2}
and multiplying the inverse factors (in reverse order) from the right
by the matrix $\Ec$, respectively $\Bc$, from~\mr{CAC}, we obtain the
relation
\[
\bea{rl}
\Mc := &\!\!\!\! \bigl(s_0\, \Ec - \Ac\bigr)^{-1} \Ec \\[12pt]
 = &\!\!\!\! \Matrix{I & 0 & \cdots & \cdots & 0 \cr
          \noalign{\vskip2pt}
             s_0 I & I & 0 & \cdots & 0\cr
             s_0^2 I & s_0 I & I & \ddots & \vdots \cr
             \vdots & \ddots &  \ddots  & \ddots  & 0 \cr
        s_0^{l-1} I & \cdots & s_0^2 I & s_0 I & I}
  \Matrix{M^{(1)} & M^{(2)} & M^{(3)} & \cdots & M^{(l)}\cr
          \noalign{\vskip2pt}
          -I & 0 & 0 & \cdots & 0 \cr
          0 & -I & 0 & \cdots & 0 \cr
     \vdots & \ddots & \ddots & \ddots & \vdots \cr
     0 & \cdots & 0 & -I  & 0},
\ea
\]
respectively 
\[
\Rc := 
\bigl(s_0\, \Ec - \Ac\bigr)^{-1} \Bc
    = \Matrix{I \cr
              s_0 I \cr \noalign{\vskip2pt}
              s_0^2 I \cr \noalign{\vskip2pt}
              \vdots \cr \noalign{\vskip2pt}
             s_0^{l-1} I} \bigl(P(s_0) \bigr)^{-1} B.
\]
The first relation is readily rewritten in the form~\mr{PIA},
and the second relation is just~\mr{PIC}.
Thus the proof is complete.

\section*{Appendix B.}

In this appendix, we establish the representations~\mr{STA} and~\mr{STM}
for the case that $P_{-1}$ is of the form~\mr{AF1} and~\mr{AF2}, respectively.

First assume that $P_{-1}$ is given by~\mr{AF1}. 
Using~\mr{AF1},~\mr{CAM}, and~\mr{QS0}, one readily verifies that
\[
\bea{rl}
s_0\, \Ec - \Ac = &\!\!\!\! 
  \Matrix{s_0 P_1 + P_0 & F_1 G \cr \noalign{\vskip4pt}
              - F_2^H &  s_0 I_{\hat{n}_0}}\\[16pt]  
= &\!\!\!\! \Matrix{Q_0 &  F_1 G \cr \noalign{\vskip4pt}
               0 &  s_0 I_{\hat{n}_0}}
     \Matrix{I_{n_0} &  0 \cr \noalign{\vskip4pt}
             - \szinv F_2^H &  I_{\hat{n}_0}}.
\ea
\]
It follows that
\beq{XAB1}{
 \bigl(s_0\, \Ec - \Ac\bigr)^{-1} =
 \Matrix{I_{n_0} &  0 \cr \noalign{\vskip4pt}
            \szinv F_2^H &  I_{\hat{n}_0}}
\Matrix{Q_0^{-1} &  -\szinv Q_0^{-1} F_1 G \cr \noalign{\vskip4pt}
               0 &  \szinv I_{\hat{n}_0}}.
}
By multiplying~\mr{XAB1} from the right by the matrix $\Ec$, respectively $\Rc$, 
from~\mr{CAM}, we obtain the relations stated in~\mr{STA}.

Next, we assume that $P_{-1}$ is given by~\mr{AF2}. 
Recall that the matrix $G$ is nonsingular.
In this case, we have
\[
\bea{rl}
s_0\, \Ec - \Ac = &\!\!\!\! 
  \Matrix{s_0 P_1 + P_0 & F_1 \cr \noalign{\vskip4pt}
              - F_2^H &  s_0 G}\\[16pt]
= &\!\!\!\! \Matrix{Q_0 &  F_1 \cr \noalign{\vskip4pt}
               0 &  s_0 G}
     \Matrix{I_{n_0} &  0 \cr \noalign{\vskip4pt}
             - \szinv G^{-1} F_2^H &  I_{\hat{n}_0}}.
\ea
\]
It follows that
\beq{XAB2}{
 \bigl(s_0\, \Ec - \Ac\bigr)^{-1} =
 \Matrix{I_{n_0} &  0 \cr \noalign{\vskip4pt}
            \szinv G^{-1} F_2^H &  I_{\hat{n}_0}}
\Matrix{Q_0^{-1} &  -\szinv Q_0^{-1} F_1 G^{-1} \cr \noalign{\vskip4pt}
               0 &  \szinv G^{-1}}.
}
By multiplying~\mr{XAB2} from the right by the matrix $\Ec$, respectively $\Rc$, 
from~\mr{CAM}, we obtain the relations stated in~\mr{STM}.

\end{document}